\documentclass[12pt]{article}
\usepackage{latexsym,amsfonts,amssymb}
\usepackage{threeparttable}
\usepackage{amsmath}
\usepackage{fancyhdr}
\usepackage{tocloft}
\usepackage{extarrows}

\usepackage[dvips]{color}
\usepackage{colordvi,multicol}

\newtheorem{theorem}{Theorem}[section]
\newtheorem{corollary}[theorem]{Corollary}

\newtheorem{proposition}[theorem]{Proposition}
\newtheorem{definition}[theorem]{Definition}
\newtheorem{remark}[theorem]{Remark}

\newcommand{\nn}{\nonumber}

\begin{document}

\title{\bf Jensen's Inequality for  Backward SDEs Driven by $G$-Brownian motion}
\author{Ze-Chun Hu\thanks{Corresponding
author: Department of Mathematics, Nanjing University, Nanjing
210093, China\vskip 0cm E-mail address: huzc@nju.edu.cn}\ \  and Zhen-Ling Wang \\
 {\small Nanjing University}}
 \date{}
 \maketitle

\noindent{\bf Abstract}\quad In this note, we consider Jensen's inequality for the nonlinear expectation associated with backward SDEs driven by $G$-Brownian motion ($G$-BSDEs for short).
At first, we give a necessary and sufficient condition for $G$-BSDEs
under which one-dimensional Jensen inequality holds. Second,  we prove that for $n>1$, the $n$-dimensional Jensen inequality holds for any nonlinear
expectation if and only if the nonlinear expectation is linear, which is essentially due to Jia (Arch. Math. 94 (2010), 489-499).
 As a consequence, we give a necessary and sufficient condition for $G$-BSDEs
under which the $n$-dimensional Jensen inequality holds.

\smallskip



\smallskip



\section{Introduction}

It's well known that  backward stochastic differential equations (BSDEs in short) play a very important role in stochastic analysis, finance and etc. We refer to a survey paper of Peng \cite{Pe10b} for more details of the theoretical studies and
applications to, e.g., stochastic controls, optimizations, games and finance.

Peng \cite{Pe04}-\cite{Pe10} defined the $G$-expectations, $G$-Brownian motions and built It\^{o}'s type stochastic calculus. As to the classic setting, it's important to study BSDEs under $G$-expectation, i.e. BSDEs  driven by $G$-Brownian motions ($G$-BSDE for short). By Hu et al. \cite{HJPS12a},  a general $G$-BSDE is to find a triple of processes $(Y,Z,K)$, where $K$ is a decreasing $G$-martingale, satisfying
\begin{eqnarray}\label{main}
Y_t&=&\xi+\int_t^Tf(s,Y_s,Z_s)ds+\int_t^Tg(s,Y_s,Z_s)d\langle B\rangle_s\nonumber\\
&&-\int_t^TZ_sdB_s-(K_T-K_t).
\end{eqnarray}

When the generator $f$ in (\ref{main}) is independent of $z$ and $g=0$, the above prolem can be equivalently formulated as
$$
Y_t=\hat{\mathbb{E}}_t[\xi+\int_t^Tf(s,Y_s)ds].
$$
The existence and uniqueness of such fully nonlinear BSDE was obtained in Peng \cite{Pe05,Pe07b,Pe10}. Soner, Touzi and Zhang \cite{STZ12} have proved the existence and uniqueness  for a type of fully nonlinear BSDE, called 2BSDE, whose generator can contain $Z$-term.

 For the general $G$-BSDE (\ref{main}), Hu et al.   proved the existence and uniqueness in \cite{HJPS12a}, and studied comparison theorem, nonlinear Feynman-Kac formula and Girsanov transformation in \cite{HJPS12b}. He and Hu \cite{HH13} obtained a representation theorem for the generators of $G$-BSDEs and used the representation theorem to get a converse comparison theorem for $G$-BSDEs and some equivalent results for the nonlinear expectations generated by $G$-BSDEs. Peng and Song \cite{PS13} introduced a new notion of $G$-expectation-weighted Sobolev spaces ($G$-Sobolev space for short), and proved that  $G$-BSDEs are in fact path dependent PDEs in the corresponding $G$-Sobolev spaces.

In this note, we study Jensen's inequality for $G$-BSDEs.  For Jensen's inequality for $g$-expectation associated classical BSDEs, we refer to  Briand et al. \cite{BHMP00}, Chen et al. \cite{CKJ03}, Jiang and Chen \cite{JC04}, Hu \cite{Hu05}, Jiang \cite{Ji06}, Fan \cite{Fa09},  Jia \cite{Ji10}, Jia and Peng \cite{JP10} and the references therein.

Recently, Guessab and Schmeisser \cite{GS13} considered the $d$-dimensional Jensen inequality
\begin{eqnarray*}
T[\psi(f_1,\cdots,f_d)]\geq \psi(T[f_1],\cdots,T[f_d]),
\end{eqnarray*}
where $T$ is a functional, $\psi$ is a convex function defined on a closed convex set $K\subset \mathbb{R}^d$, and $f_1,\cdots,f_d$ are from some linear space of functions.
Among other things, the authors showed that if we exclude three types of convex sets $K$, then Jensen's inequality holds for a sublinear functional $T$ if and only if $T$ is linear, positive, and satisfies $T[1]=1$, i.e. $T$ is a linear expectation.

The rest of this note is organized as follows. In Section 2, we give some preliminaries about $G$-expectation and $G$-BSDEs.
In Section 3,  we consider Jensen's inequality for the nonlinear expectation driven by $G$-BSDEs.
In Subsection 3.1, we follow the method of Hu \cite{Hu05} and apply the comparision theorem, the converse comparison theorem in He and Hu \cite{HH13} to give a necessary and sufficient condition for $G$-BSDEs
under which one-dimensional Jensen inequality holds. In Subsection 3.2,  we prove that for $n>1$, the $n$-dimensional Jensen inequality holds for any nonlinear
expectation if and only if the nonlinear expectation is linear, which is essentially due to Jia \cite{Ji10}, and as a consequence, we give a necessary and sufficient condition for $G$-BSDEs
under which the $n$-dimensional Jensen inequality holds.

\section{Preliminaries}
In this section, we review some basic notions and results of $G$-expectation, the related spaces of random variables, and  $G$-BSDE. The readers may refer to \cite{Pe10}, \cite{HJPS12a} and  \cite{HJPS12b}  for more details.

\begin{definition}\label{def2.1} Let $\Omega$ be a given set and let $\mathcal{H}$ be a linear space of real valued function defined on $\Omega$, and satisfy:
$(i)$ for each constant $c$,  $c\in \mathcal{H}$; $(ii)$ if $X\in\mathcal{H}$, then $|X|\in\mathcal{H}$. The space $\mathcal{H}$ can be considered as the space of random variables.
A sublinear expectation $\hat{\mathbb{E}}$ is a functional $\hat{\mathbb{E}}:\mathcal{H}\rightarrow\mathbb{R}$ satisfying\\
$(i)$ Monotonicity:~$\hat{\mathbb{E}}[X]\ge\hat{\mathbb{E}}[Y]$, if $X\ge Y;$\\
$(ii)$ Constant preserving:~$\hat{\mathbb{E}}[c]=c$, for $c\in\mathbb{R};$\\
$(iii)$ Sub-additivity:~$\hat{\mathbb{E}}[X+Y]\le\hat{\mathbb{E}}[X]+\hat{\mathbb{E}}[Y]$, for each $X,~Y\in\mathcal{H};$\\
$(iv)$ Positive homogeneity:~$\hat{\mathbb{E}}[\lambda X]=\lambda\hat{\mathbb{E}}[X]$, for $\lambda\geq 0.$\\
The triple $(\Omega,\mathcal{H},\hat{\mathbb{E}})$ is called a sublinear expectation space. If $(i)$ and $(ii)$ are satisfied, $\hat{\mathbb{E}}$ is called a nonlinear expectation  and the triple $(\Omega,\mathcal{H},\hat{\mathbb{E}})$ is called a nonlinear expectation space.
\end{definition}

\begin{definition}
Let $X_{1}$ and $X_{2}$ be two n-dimensional random vectors defined in sublinear expectation spaces $(\Omega,\mathcal{H},\hat{\mathbb{E}}_{1})$ and $(\Omega,\mathcal{H},\hat{\mathbb{E}}_{2})$ respectively. They are called identically distributed, denoted by $X_{1}\overset{d}{=}X_{2}$, if
$
\hat{\mathbb{E}}_{1}[\varphi(X_{1})]=\hat{\mathbb{E}}_{2}[\varphi(X_{2})],
$
for all $\varphi \in C_{b.Lip}(\mathbb{R}^{n})$, where $C_{b.Lip}(\mathbb{R}^{n})$ denotes the space of all bounded and Lipschitz functions on $\mathbb{R}^n$.
\end{definition}

\begin{definition}
In a sublinear expectation space $(\Omega,\mathcal{H},\hat{\mathbb{E}})$, a random vector $Y\in{\mathcal{H}^n}$ is said to be independent of another random vector $X\in{\mathcal{H}^m}$ under $\hat{\mathbb{E}}[\cdot]$, denoted by $Y\perp X$, if for all $\varphi \in C_{b.Lip}{(\mathbb{R}^{n+m})}$ one has
$
\hat{\mathbb{E}}[\varphi(X,Y)]=\hat{\mathbb{E}}[\hat{\mathbb{E}}[\varphi(x,Y)]|_{x=X}].
$
\end{definition}

\begin{definition}
$(G$-normal distribution$)$
A d-dimensional random vector $X=(X_{1},\cdots,X_{d})$ in sublinear expectation space $(\Omega,\mathcal{H},\hat{\mathbb{E}})$ is called $G$-normally distributed if for each $a,~b\ge 0$, one has
$
aX+b\bar{X}\overset{d}{=}\sqrt{a^2+b^2}X,
$
where $\bar{X}$ is an independent copy of $X$, i.e. $\bar{X}\overset{d}{=}X$ and $\bar{X}\perp X.$ Here, the  letter $G$  denotes the function
\begin{align}
G(A):=\hat{\mathbb{E}}[\frac{1}{2}\langle AX,X\rangle]:\mathbb{S}_{d}\rightarrow\mathbb{R},\nn
\end{align}
where $\mathbb{S}_{d}=\{A|A\mbox{ is $d\times d$ symmetric matrix}\}.$
\end{definition}

Peng \cite{Pe08b} proved that  $X=(X_1,\cdots,X_d)$ is $G$-normally distributed if and only if for each $\varphi \in C_{b.Lip}{(\mathbb{R}^d)},~u(t,x):=\hat{\mathbb{E}}[\varphi(x+\sqrt{t}X)],~(t,x)\in [0,\infty)\times \mathbb{R}^d$, is the solution of the following $G$-heat equation:
\begin{align}
\partial _t u-\textit{G}(D^2_{x}u)=0,~~~u(0,x)=\varphi.\nn
\end{align}

The function $G(\cdot):\mathbb{S}_{d}\rightarrow\mathbb{R}$ is a monotonic, sublinear mapping on $\mathbb{S}_{d}$ and
$
G(A):=\hat{\mathbb{E}}[\frac{1}{2}\langle AX,X\rangle]\leq\frac{1}{2}|A|\hat{\mathbb{E}}[|X|^{2}],
$
which implies that there exists a bounded, convex, and closed subset $\Gamma\subset\mathbb{S}^{+}_{d}$ such that
\begin{align}
G(A)=\frac{1}{2}\sup_{\gamma\in \Gamma}{\rm tr}[\gamma A],\nn
\end{align}
where  $\mathbb{S}^{+}_{d}$ denotes the collection of nonnegative elements in $\mathbb{S}_{d}$. In this note, we only consider nondegenerate $G$-normal distribution; that is, there exists some $\sigma^{2}>0$ such that
$
G(A)-G(B)\geq \sigma^{2}{\rm tr}[A-B]$ for any $A\geq B.$

\begin{definition}
 (i) Let $\Omega=C_{0}^d(\mathbb{R}^{+})$ denote the space of $\mathbb{R}^{d}$-valued continuous functions on $[0,\infty)$ with $\omega_{0}=0$ and $B_{t}(\omega)=\omega_{t}$ be the canonical process. For each fixed $T\in[0,\infty)$, we set
\begin{align}
L_{ip}(\Omega_{T}):=\{\varphi(B_{t_{1}\wedge T},\cdots,B_{t_{n}\wedge T}):n\in \mathbb{N},~t_{1},\cdots,t_{n}\in[0,\infty),~\varphi\in C_{b.Lip}(\mathbb{R}^{d\times n})\}.\nn
\end{align}
It is clear that $L_{ip}(\Omega_{t})\subseteq L_{ip}(\Omega_{T})$ for $t\leq T$. We also set
$
L_{ip}(\Omega):=\bigcup_{n=1}^{\infty}L_{ip}(\Omega_{n}).
$
Let $G:\mathbb{S}_{d}\rightarrow \mathbb{R}$ be a given monotonic and sublinear function. $G$-expectation is a sublinear expectation defined by
\begin{align}
\hat{\mathbb{E}}[X]=\bar{\mathbb{E}}[\varphi(\sqrt{t_{1}-t_{0}}\xi_{1},\cdots,\sqrt{t_{m}-t_{m-1}}\xi_{m})]£¬\nn
\end{align}
for all $X\in L_{ip}(\Omega)$ with $X=\varphi(B_{t_{1}}-B_{t_{0}},B_{t_{2}}-B_{t_{1}},\cdots,B_{t_{m}}-B_{t_{m-1}})$, where $\xi_{1},\cdots,\xi_{m}$ is identically distributed $d$-dimensional $G$-normally distributed random vectors in a sublinear expectation space $(\bar{\Omega},\bar{\mathcal{H}},\bar{\mathbb{E}})$ such that $\xi_{i+1}$ is independent of $(\xi_{1},\cdots,\xi_{i})$ for every $i=1,\cdots,m-1$. The corresponding canonical process $B_{t}=(B_{t}^{i})_{i=1}^{d}$ is called a $G$-Brownian motion.

 (ii) For each fixed $t\in[0,\infty)$, the conditional $G$-expectation $\hat{\mathbb{E}}_{t}[\cdot]$ for $X=\varphi(B_{t_{1}}-B_{t_{0}},B_{t_{2}}-B_{t_{1}},\cdots,B_{t_{m}}-B_{t_{m-1}})\in L_{ip}(\Omega)$, where without loss of generality we suppose $t=t_{i}$, $1\leq i\leq m$, is defined by
\begin{align}
\hat{\mathbb{E}}_{t}[\varphi(B_{t_{1}}-B_{t_{0}},B_{t_{2}}-B_{t_{1}},\cdots,B_{t_{m}}-B_{t_{m-1}})]
=\psi(B_{t_{1}}-B_{t_{0}},B_{t_{2}}-B_{t_{1}},\cdots,B_{t_{i}}-B_{t_{i-1}}),\nn
\end{align}
where $\psi(x_{1},\cdots,x_{i})=\hat{\mathbb{E}}[\varphi(x_{1},\cdots,x_{i},B_{t_{i+1}}-B_{t_{i}},\cdots,B_{t_{m}}-B_{t_{m-1}})].$
\end{definition}

We denote by $L_{G}^p(\Omega)$, $p\geq1$, the completion of $L_{ip}(\Omega)$ under the norm
$
\|X\|_{p,G}=(\hat{\mathbb{E}}[|X|^{p}])^{1/p}.$
Similarly, we can define $L_{G}^p(\Omega_{T})$. It is clear that $L_{G}^q(\Omega)\subset L_{G}^p(\Omega)$ for $1\leq p\leq q$ and $\hat{\mathbb{E}}[\cdot]$ can be extended continuously to $L_{G}^1(\Omega).$

For each fixed $\mathbf{a}=(a_{1},\cdots,a_{d})\in\mathbb{R}^{d}$, $B_{t}^\mathbf{a}=\langle\mathbf{a},B_{t}\rangle$ is a $1$-dimensional $G_{\mathbf{a}}$-Brownian motion on $(\Omega,\mathcal{H},\hat{\mathbb{E}})$, where $G_{\mathbf{a}}(\alpha)=\frac{1}{2}(\sigma_{\mathbf{a}\mathbf{a}^{T}}^{2}\alpha^{+}-\sigma_{-\mathbf{a}\mathbf{a}^{T}}^{2}\alpha^{-}),~\sigma_{\mathbf{a}\mathbf{a}^{T}}^{2}=2G(\mathbf{a}\mathbf{a}^{T})=\hat{\mathbb{E}}[\langle\mathbf{a},B_{1}\rangle^{2}],~\sigma_{-\mathbf{a}\mathbf{a}^{T}}^{2}=-2G(-\mathbf{a}\mathbf{a}^{T})=-\hat{\mathbb{E}}[-\langle\mathbf{a},B_{1}\rangle^{2}]$. In particular, for each $t,~s\geq 0,~B_{t+s}^{\mathbf{a}}-B_{t}^{\mathbf{a}}\overset{d}{=}N({0}\times[s\sigma_{-\mathbf{a}\mathbf{a}^{T}}^{2},s\sigma_{\mathbf{a}\mathbf{a}^{T}}^{2}]).$

Let $\pi_{T}^{N}=\{t_{0}^{N},t_{1}^{N},\cdots,t_{N}^{N}\},~N=1,2,\cdots$, be a sequence of partitions of $[0,t]$ such that $\mu(\pi_{T}^N)=\max\{|t_{i+1}-t_{i}|:~i=0,1,\cdots,N-1\}\rightarrow 0.$
The quadratic variation process of $\langle B^{\mathbf{a}}\rangle$ is defined by
$$\langle B^{\mathbf{a}}\rangle_{t}:=\lim_{\mu(\pi_{T}^{N})\rightarrow 0}\sum_{k=0}^{N-1}(B_{t_{k+1}^{N}}^{\mathbf{a}}-B_{t_{k}^{N}}^{\mathbf{a}})^{2}=(B_{t}^{\mathbf{a}})^{2}-2\int_{0}^{t}B_{s}^{\mathbf{a}}dB_{s}^{\mathbf{a}}.$$
For each fixed $\mathbf{a},\mathbf{\bar{a}}\in\mathbb{R}^{d}$, the mutual variation process of $B^{\mathbf{a}}$ and $B^{\mathbf{\bar{a}}}$ is defined by
\begin{align}
\langle B^{\mathbf{a}},B^{\mathbf{\bar{a}}}\rangle_{t}:=\frac{1}{4}[\langle B^{\mathbf{a}}+B^{\mathbf{\bar{a}}}\rangle_{t}-\langle B^{\mathbf{a}}-B^{\mathbf{\bar{a}}}\rangle_{t}]=\frac{1}{4}[\langle B^{\mathbf{a}+\bar{\mathbf{a}}}\rangle_{t}-\langle B^{\mathbf{a}-\bar{\mathbf{a}}}\rangle_{t}].\nn
\end{align}

\begin{definition}
For fixed $T\geq 0$, let $ M_{G}^{0}(0,T)$ be the collection of process in the following form: for a given partition
$\pi_{T}=\{t_{0},t_{1},\cdots,t_{N}\}$ of $[0,T]$,
\begin{align}
\eta_{t}(\omega)=\sum_{k=0}^{N-1}\xi_{k}(\omega)\mathbf{I}_{[t_{k},t_{k+1})}(t),\nn
\end{align}
where $\xi_{k}\in L_{G}^p(\Omega_{t_{k}}),~k=0,1,2,\cdots,N-1$. For $p\geq 1$, we denote by $H_{G}^{p}(0,T),~M_{G}^{p}(0,T)$ the completion of $M_{G}^{0}(0,T)$ under the norms
$
\|\eta\|_{H_{G}^{p}}=\{\hat{\mathbb{E}}[(\int_{0}^{T}|\eta_{t}|^{2}dt)^{p/2}]\}^{1/p},
\|\eta\|_{M_{G}^{p}}=\{\hat{\mathbb{E}}[\int_{0}^{T}|\eta_{t}|^{p}dt]\}^{1/p},
$ respectively.
\end{definition}

Let $S_{G}^{0}(0,T)=\{h(t,B_{t_{1}\wedge t},\cdots,B_{t_{n}\wedge t}):t_{1},\cdots,t_{n}\in[0,T],~h\in C_{b.Lip}(\mathbb{R}^{n+1})\}$. For $p\geq1$, denote by $S_{G}^{p}(0,T)$ the completion of $S_{G}^{0}(0,T)$ under the norm
$
\|\eta\|_{S_{G}^{p}}=\{\hat{\mathbb{E}}[\sup_{t\in[0,T]}|\eta_{t}|^{p}]\}^{\frac{1}{p}}.
$

We consider the following type of $G$-BSDEs (in this note we always use Einstein convention):
\begin{eqnarray}\label{eq3}
Y_{t}&=&\xi+\int_{t}^{T}f(s,Y_{s},Z_{s})ds+\int_{t}^{T}g_{ij}(s,Y_{s},Z_{s})d\langle B^{i},B^{j}\rangle_{s}\nn\\
&&-\int_{t}^{T}Z_{s}dB_{s}-(K_{T}-K_{t}),
\end{eqnarray}
where
$$
f(t,\omega,y,z),~g_{ij}(t,\omega,y,z):[0,T]\times\Omega_{T}\times\mathbb{R}\times\mathbb{R}^{d}\rightarrow\mathbb{R},
$$
satisfy the following properties:
\begin{itemize}
\item[(H1)] There exists some $\beta>1$ such that for any $y,~z,~f(\cdot,\cdot,y,z),~g_{ij}(\cdot,\cdot,y,z)\in M_{G}^{\beta}(0,T);$

\item[(H2)] There exists some $L>0$ such that
\begin{align}
|f(t,y,z)-f(t,y',z')|+\sum_{i,j=1}^{d}|g_{ij}(t,y,z)-g_{ij}(t,y',z')|\leq L(|y-y'|+|z-z'|).\nn
\end{align}
\end{itemize}

Denote by $\mathfrak{S}_{G}^{\alpha}(0,T)$ the completion of processes $(Y,Z,K)$ such that $Y\in S_{G}^{\alpha}(0,T),~Z\in H_{G}^{\alpha}(0,T;\mathbb{R}^{d}),~K$ is a decreasing $G$-martingale with $K_{0}=0$ and $K_{T}\in L_{G}^{\alpha}(\Omega_{T}).$

\begin{definition}
Let $\xi\in L_{G}^{\beta}(\Omega_{T})$ and $f$ and $g_{ij}$ satisfy (H1) and (H2) for some $\beta>1$. A triplet of processes $(Y,Z,K)$ is called a solution of \eqref{eq3} if for some $1<\alpha\leq\beta$ the following properties hold:
\begin{itemize}
\item[(a)] $(Y,Z,K)\in\mathfrak{S}_{G}^{\alpha}(0,T);$
\item[(b)] $Y_{t}=\xi+\int_{t}^{T}f(s,Y_{s},Z_{s})ds+\int_{t}^{T}g_{ij}(s,Y_{s},Z_{s})d\langle B^{i},B^{j}\rangle_{s}-\int_{t}^{T}Z_{s}dB_{s}-(K_{T}-K_{t}).$
\end{itemize}
\end{definition}

\begin{theorem}\label{Th1} (\cite{HJPS12a})
Assume that $\xi\in L_{G}^{\beta}(\Omega_{T})$ and $f$ and $g_{ij}$ satisfy $(H1)$ and $(H2)$ for some $\beta>1$. Then, equation \eqref{eq3} has a unique solution $(Y,Z,K)$. Moreover, for any $1<\alpha<\beta$, one has $Y\in S_{G}^{\alpha}(0,T),~Z\in H_{G}^{\alpha}(0,T;\mathbb{R}^{d})$ and $K_{T}\in L_{G}^{\alpha}(\Omega_{T}).$
\end{theorem}

In this note, we also need the following assumptions for $G$-BSDE \eqref{eq3} (see He and Hu \cite{HH13}).

\begin{itemize}
\item[(H3)] For each fixed $(\omega,y,z)\in \Omega_{T}\times\mathbb{R}\times\mathbb{R}^{d},~t\rightarrow f(t,\omega,y,z)$ and $t\rightarrow g_{ij}(t,\omega,y,z)$ are continuous.

\item[(H4)] For each fixed $(t,y,z)\in [0,T)\times\mathbb{R}\times\mathbb{R}^{d},~f(t,y,z),~g_{ij}(t,y,z)\in L_{G}^{\beta}(\Omega_{t})$, and
\begin{align}
\lim_{\varepsilon\rightarrow 0+}\frac{1}{\varepsilon}\hat{\mathbb{E}}\left[\int_{t}^{t+\varepsilon}\left(|f(u,y,z)-f(t,y,z)|^{\beta}+\sum_{i,j=1}^{d}|g_{ij}(u,y,z)-g_{ij}(t,y,z)|^{\beta}\right)du\right]=0.
\end{align}

\item[(H5)] For each fixed $(t,\omega,y)\in [0,T]\times\Omega_{T}\times\mathbb{R},~f(t,\omega,y,0)=g_{ij}(t,\omega,y,0)=0.$

\end{itemize}


\section{Jensen's inequality for   $G$-BSDEs}\setcounter{equation}{0}

We consider the following  $G$-BSDE:
\begin{align}\label{eq15}
Y_{t}=&\xi+\int_{t}^{T}f(s,Y_{s},Z_{s})ds+\int_{t}^{T}g_{ij}(s,Y_{s},Z_{s})d\langle B^{i},B^{j}\rangle_{s}\nn\\
&-\int_{t}^{T}Z_{s}dB_{s}-(K_{T}-K_{t}),
\end{align}where $g_{ij}=g_{ji}$, and $f$ and $g_{ij}$ satisfy the conditions (H1)-(H5). Define $\tilde{\mathbb{E}}_{t}[\xi]=Y_{t}.$

\subsection{One-dimensional Jensen inequality}

\begin{theorem}\label{thm}
The following two statements are equivalent:\\
(i) Jensen's inequality holds, i.e, for each $\xi\in L_{G}^{2}(\Omega_{T})$, and any convex function $h:\mathbb{R}\rightarrow\mathbb{R}$, if $h(\xi)\in L_{G}^{2}(\Omega_{T})$, then
\begin{align}\label{thm-a}
\tilde{\mathbb{E}}_{t}[h(\xi)]\geq h(\tilde{\mathbb{E}}_{t}[\xi]),~~~\forall t\in[0,T].
\end{align}
(ii)  $\forall \lambda,~\mu\in\mathbb{R},~\lambda\neq0,~\forall (t,y,z)\in [0,T]\times\mathbb{R}\times\mathbb{R}^{d}$,
\begin{align}\label{thm-b}
\lambda f(t,y,z)-f(t,\lambda y+\mu,\lambda z)+2G((\lambda g_{ij}(t,y,z)-g_{ij}(t,\lambda y+\mu,\lambda z))_{i,j=1}^{d})\leq 0,~~~q.s.
\end{align}
\end{theorem}

\noindent{\bf Proof.} The idea of the proof comes from Theorem 3.1 of \cite{Hu05}.

$(i)\Rightarrow (ii):$ For fixed $\lambda\neq0$ and $\mu$, we define a convex function $h(x)=\lambda x+\mu$. Let  $(Y_{t},Z_{t},K_{t})$ be the unique solution of the $G$-BSDE \eqref{eq15}. Define $Y_{t}'=\lambda Y_{t}+\mu,~Z_{t}'=\lambda Z_{t},~K_{t}'=\lambda K_{t}$. Then $(Y_{t}',Z_{t}',K_{t}')$ is the unique solution of the following $G$-BSDE:
\begin{align}\label{thm-c}
Y_{t}'=&h(\xi)+\int_{t}^{T}f'(s,Y_{s}',Z_{s}')ds+\int_{t}^{T}g_{ij}'(s,Y_{s}',Z_{s}')d\langle B^{i},B^{j}\rangle_{s}\nn\\
&-\int_{t}^{T}Z_{s}'dB_{s}-(K_{T}'-K_{t}'),
\end{align}
where $f'(t,y,z)=\lambda f(t,\frac{y-\mu}{\lambda},\frac{z}{\lambda}),~g_{ij}'(t,y,z)=\lambda g_{ij}(t,\frac{y-\mu}{\lambda},\frac{z}{\lambda})$.

Denote $\tilde{\mathbb{E}}_{t}'[h(\xi)]=Y_{t}'$. By \eqref{thm-a}, we get
\begin{align}\label{thm-d}
\tilde{\mathbb{E}}_{t}[h(\xi)]\geq h(\tilde{\mathbb{E}}_{t}[\xi])=\lambda Y_{t}+\mu=Y_{t}'=\tilde{\mathbb{E}}_{t}'[h(\xi)].
\end{align}
For any $\eta\in L_{G}^{2}(\Omega_{T})$, put $\xi=h^{-1}(\eta)$. Then we have by (\ref{thm-d})
\begin{align}
\tilde{\mathbb{E}}_{t}[\eta]\geq\tilde{\mathbb{E}}_{t}'[\eta].\nn
\end{align}
By  the converse comparison theorem \cite[Theorem 15]{HH13}, we obtain that
\begin{align}
(f'-f)(t,y',z')+2G((g_{ij}'-g_{ij})_{i,j=1}^{d})(t,y',z')\leq 0\ q.s.,\nn
\end{align}
which implies
\begin{align}
&f'(t,y',z')-f(t,y',z')+2G((g_{ij}'(t,y',z')-g_{ij}(t,y',z'))_{i,j=1}^{d})\nn\\
&=\lambda f(t,\frac{y'-\mu}{\lambda},\frac{z'}{\lambda})-f(t,y',z')+2G((\lambda g_{ij}(t,\frac{y'-\mu}{\lambda},\frac{z'}{\lambda})-g_{ij}(t,y',z'))_{i,j=1}^{d})\nn\\
&\xlongequal[z:=\frac{z'}{\lambda}]{y:=\frac{y'-\mu}{\lambda}}\lambda f(t,y,z)-f(t,\lambda y+\mu,\lambda z)+2G((\lambda g_{ij}(t,y,z)-g_{ij}(t,\lambda y+\mu,\lambda z))_{i,j=1}^{d})\nn\\
&\leq 0,~~~q.s.\nn
\end{align}
Hence $(ii)$ holds.

$(ii) \Rightarrow (i):$ First,  take a linear function $h(x)=\lambda x+\mu$ where $\lambda\neq0$. Let $(Y_{t},Z_{t},K_{t})$ be the unique solution of $G$-BSDE \eqref{eq15}, and denote $Y_{t}'=\lambda Y_{t}+\mu,~Z_{t}'=\lambda Z_{t},~K_{t}'=\lambda K_{t}$. Then $(Y_{t}',Z_{t}',K_{t}')$ is the unique solution of $G$-BSDE \eqref{thm-c}. Let $f',g_{ij}'$ be defined as in \eqref{thm-c}. Then by $(ii)$, we have
$$(f'-f)(t,y,z)+2G((g_{ij}'-g_{ij})_{i,j=1}^{d})(t,y,z)\leq 0\ q.s.,$$
which together with  the comparision theorem \cite[Proposition 13]{HH13}  implies that
\begin{align}\label{thm-e}
\tilde{\mathbb{E}}_{t}[h(\xi)]\geq \tilde{\mathbb{E}}_{t}'[h(\xi)]=Y_{t}'=\lambda Y_{t}+\mu=\lambda\tilde{\mathbb{E}}_{t}[\xi]+\mu=h(\tilde{\mathbb{E}}_{t}[\xi]).
\end{align}

For any convex function $h$, there exists a countable set $D$ in $\mathbb{R}^{2}$, such that
\begin{align}\label{thm-f}
h(x)=\sup_{(\lambda,\mu)\in D}(\lambda x+\mu).
\end{align}
By (\ref{thm-e}) and (\ref{thm-f}), we have
\begin{align}
\tilde{\mathbb{E}}_{t}[h(\xi)]= \tilde{\mathbb{E}}_{t}[\sup_{(\lambda,\mu)\in D}(\lambda \xi+\mu)]\geq\sup_{(\lambda,\mu)\in D}(\lambda\tilde{\mathbb{E}}_t[\xi]+\mu)=h(\tilde{\mathbb{E}}_{t}[\xi]),\nn
\end{align}
i.e. $(i)$ holds.
\hfill\fbox

\begin{remark}
 (i) If $f$ and $g_{ij}$ are independent of $y$, then the  condition of \eqref{thm-b}  becomes
\begin{align}
\lambda f(t,z)-f(t,\lambda z)+2G((\lambda g_{ij}(t,z)-g_{ij}(t,\lambda z))_{i,j=1}^{d})\leq 0,\ q.s.\nn
\end{align}
(ii) If $g_{ij}\equiv0$, then the  condition of \eqref{thm-b}  becomes
\begin{align}\label{rem}
f(t,\lambda y+\mu,\lambda z)\geq \lambda f(t,y,z),\ q.s.
\end{align}
Taking $\lambda=1$, then $f(t,y+\mu, z)\geq  f(t,y,z),\ q.s.$, which implies that $f$ is independent of $y$. Thus \eqref{rem} becomes $f(t,\lambda z)\geq \lambda f(t,z),\ q.s.$ This is just the condition in Hu \cite[Theorem 3.1]{Hu05}.
\end{remark}

\subsection{Multi-dimensional Jensen inequality}

At first, we prove a  result for any nonlinear expectation, which is essentially due to  Jia (see \cite[Theorem 3.3]{Ji10}).

\begin{theorem}\label{thm3.3}
Assume that $n>1$ and $(\Omega,\mathcal{H},\hat{\mathbb{E}})$ is a nonlinear expectation space defined by Definition \ref{def2.1}. Then the following two claims are equivalent:\\
(a)  $\hat{\mathbb{E}}$ is linear, i.e., for any $\lambda,\gamma\in \mathbb{R},X,Y\in \mathcal{H}$,
\begin{eqnarray}\label{lem3.3-a}
\hat{\mathbb{E}}[\lambda X+\gamma Y]=\lambda\hat{\mathbb{E}}[X]+\gamma \hat{\mathbb{E}}[Y];
\end{eqnarray}
(b) the $n$-dimensional Jensen inequality for nonlinear expectation $\hat{\mathbb{E}}$ holds, i.e. for each $X_i\in \mathcal{H}(i=1,\cdots,n)$ and convex function $h:\mathbb{R}^n\to \mathbb{R}$, if $h(X_1,\cdots,X_n)\in \mathcal{H}$, then
$$
\hat{\mathbb{E}}[h(X_1,\cdots,X_n)]\geq h(\hat{\mathbb{E}}[X_1],\cdots,\hat{\mathbb{E}}[X_n]).
$$
\end{theorem}
{\bf Proof.} The proof of \cite[Theorem 3.3]{Ji10} can be moved to this case. For the reader's convenience, we spell out the details.

$(b)\Rightarrow (a)$: For any $(\lambda_1,\cdots,\lambda_n)\in \mathbb{R}^n$, by (b) we have that
\begin{eqnarray}\label{lem3.3-b}
\hat{\mathbb{E}}\left[\sum_{i=1}^n\lambda_iX_i\right]\geq \sum_{i=1}^N\lambda_i\hat{\mathbb{E}}[X_i].
\end{eqnarray}
Taking $\lambda_1> 0,\lambda_j=0,j=2,\cdots,n$, we get that
$$
\hat{\mathbb{E}}\left[\lambda_1X_1\right]\geq \lambda_1\hat{\mathbb{E}}[X_1]\geq \lambda_1\cdot \frac{1}{\lambda}\hat{\mathbb{E}}[\lambda X_1]=\hat{\mathbb{E}}\left[\lambda_1X_1\right],
$$
which together with $\hat{\mathbb{E}}[0]=0$ (by (ii) in Definition \ref{def2.1}) implies that  $\hat{\mathbb{E}}$ is positively homogeneous. Put $\lambda_1=1,\lambda_2=-1$ and $\lambda_1=\lambda_2=1$
respectively, and put $\lambda_j=0$ for $j>2$ in (\ref{lem3.3-b}), we get
\begin{eqnarray*}
&&\hat{\mathbb{E}}[X_1-X_2]\geq \hat{\mathbb{E}}[X_1]-\hat{\mathbb{E}}[X_2],\ \ \hat{\mathbb{E}}[X_1+X_2]\geq \hat{\mathbb{E}}[X_1]+\hat{\mathbb{E}}[X_2].
\end{eqnarray*}
It follows that $\hat{\mathbb{E}}[X_1]\leq
\hat{\mathbb{E}}[X_2]+\hat{\mathbb{E}}[X_1-X_2]\leq \hat{\mathbb{E}}[X_2+(X_1-X_2)]=\hat{\mathbb{E}}[X_1]$. Thus we have
$\hat{\mathbb{E}}[X_1-X_2]=\hat{\mathbb{E}}[X_1]-\hat{\mathbb{E}}[X_2]$
and $\hat{\mathbb{E}}[X_1+X_2]=\hat{\mathbb{E}}[(X_1+X_2)-X_2]+\hat{\mathbb{E}}[X_2]=\hat{\mathbb{E}}[X_1]+\hat{\mathbb{E}}[X_2]$.
Hence $\hat{\mathbb{E}}$ is homogeneous and thus it's linear.

$(a)\Rightarrow (b)$: For any $(\lambda_1,\cdots,\lambda_n,\mu)\in \mathbb{R}^{n+1}$, by (a) and (ii) in Definition \ref{def2.1}, we have
\begin{eqnarray}\label{lem3.3-c}
\hat{\mathbb{E}}\left[\sum_{i=1}^n\lambda_iX_i+\mu\right]=\hat{\mathbb{E}}\left[\sum_{i=1}^n\lambda_iX_i\right]+\mu
=\sum_{i=1}^n\lambda_i\hat{\mathbb{E}}[X_i]+\mu.
\end{eqnarray}
For any convex function $h:\mathbb{R}^n\to \mathbb{R}$, there exists a countable set $D\subset \mathbb{R}^{n+1}$ such that
\begin{eqnarray}\label{lem3.3-d}
h(x)=\sup_{(\lambda_1,\cdots,\lambda_n,\mu)\in D}\left(\sum_{i=1}^n\lambda_ix_i+\mu\right).
\end{eqnarray}
By (\ref{lem3.3-c}) and (i) in Definition \ref{def2.1}, for any $(\lambda_1,\cdots,\lambda_n,\mu)\in D$, we have
$$
\hat{\mathbb{E}}[h(X_1,\cdots,X_n)]\geq \hat{\mathbb{E}}\left[\sum_{i=1}^n\lambda_iX_i+\mu\right]=\sum_{i=1}^n\lambda_i\hat{\mathbb{E}}[X_i]+\mu,
$$
which together with (\ref{lem3.3-d}) implies $(b)$. \hfill\fbox

\begin{proposition}\label{pro3.4}
Assume that $n>1$ and $t\in [0,T]$.  Then the following two claims are equivalent:\\
(i)  $\tilde{\mathbb{E}}_t$ is linear, i.e., for any $\lambda,\gamma\in \mathbb{R},X,Y\in \mathcal{H}$,
\begin{eqnarray}\label{lem3.3-a}
\tilde{\mathbb{E}}_t[\lambda X+\gamma Y]=\lambda\tilde{\mathbb{E}}_t[X]+\gamma \tilde{\mathbb{E}}_t[Y];
\end{eqnarray}
(ii) the $n$-dimensional Jensen inequality for  $\tilde{\mathbb{E}}_t$ holds, i.e. for each $X_i\in \mathcal{H}(i=1,\cdots,n)$ and convex function $h:\mathbb{R}^n\to \mathbb{R}$, if $h(X_1,\cdots,X_n)\in \mathcal{H}$, then
$$
\tilde{\mathbb{E}}_t[h(X_1,\cdots,X_n)]\geq h(\tilde{\mathbb{E}}_t[X_1],\cdots,\tilde{\mathbb{E}}_t[X_n]).
$$
\end{proposition}
{\bf Proof.} By \cite[Theorem 5.1  (1)(2)]{HJPS12b}, we know that $\tilde{\mathbb{E}}_t$ satisfies monotonicity and constant preserving.
Then all the proof of the above theorem can be moved to this case.
\hfill\fbox

\begin{corollary}
Assume that $n>1$.  Then the following two claims are equivalent:\\
(i) for any $t\in [0,T]$, the $n$-dimensional Jensen inequality for  $\tilde{\mathbb{E}}_t$ holds, i.e. for each $X_i\in \mathcal{H}(i=1,\cdots,n)$ and convex function $h:\mathbb{R}^n\to \mathbb{R}$, if $h(X_1,\cdots,X_n)\in \mathcal{H}$, then
$$
\tilde{\mathbb{E}}_t[h(X_1,\cdots,X_n)]\geq h(\tilde{\mathbb{E}}_t[X_1],\cdots,\tilde{\mathbb{E}}_t[X_n]);
$$
(ii) for any $t\in [0,T],y,y'\in \mathbb{R},z,z'\in \mathbb{R}^d,\lambda\geq 0$,
\begin{eqnarray*}
&&f(t,y+y',z+z')-f(t,y,z)-f(t,y',z')\\
&&=-2G\left((g_{ij}(t,y+y',z+z')-g_{ij}(t,y,z)-g_{ij}(t,y',z'))_{i,j=1}^d\right),
\end{eqnarray*}
and
\begin{eqnarray*}
f(t,\lambda y,\lambda z)-\lambda f(t,y,z)&=&2G\left((\lambda g_{ij}(t,y,z)-g_{ij}(t,\lambda y,\lambda z))_{i,j=1}^d\right)\\
&=&-2G\left((g_{ij}(t,\lambda y,\lambda z)-\lambda g_{ij}(t,y,z))_{i,j=1}^d\right).
\end{eqnarray*}
\end{corollary}
{\bf Proof.} By Proposition \ref{pro3.4}, we know that ($i)$ holds if and only if for any $t\in [0,T]$, $\tilde{\mathbb{E}}_t$ is linear. Then by   \cite[Proposition 17 (2)(4)]{HH13}, we obtain that $(i)$ and $(ii)$ are equivalent.\hfill\fbox

\bigskip

 { \noindent {\bf\large Acknowledgments} \vskip 0.1cm  \noindent   We are grateful to the support of
NNSFC and Jiangsu Province Basic Research Program (Natural Science
Foundation) (Grant No. BK2012720).}

\end{document}